\documentclass{article}

\usepackage[english]{babel}

\usepackage[letterpaper,top=2cm,bottom=2cm,left=3cm,right=3cm,marginparwidth=1.75cm]{geometry}
\usepackage{amsmath}
\usepackage{color}
\usepackage{graphicx}
\usepackage{placeins}
\usepackage[colorlinks=true, allcolors=black]{hyperref}
\usepackage{ulem}
\usepackage{titling}

\usepackage[affil-it]{authblk}     
\title{Seasonal Forcing in Rock–Paper–Scissors Population Dynamics}
\author[1]{Idan Sorin}
\author[2]{Michael Zaks}
\author[1]{Alexander Nepomnyashchy}

\affil[1]{Department of Mathematics, Technion–Israel Institute of Technology, Haifa 32000, Israel}
\affil[2]{Institute of Physics, Humboldt University of Berlin, Berlin, Germany}

\begin{document}
\maketitle
\begin{abstract}
We study a class of cyclic dominance models with seasonal forcing, extending the classical May–Leonard competition framework. By introducing time-periodic coefficients into the growth rates and the interaction terms, we explore how environmental seasonality influences the dynamics of three-species systems. Through analytical estimates and numerical simulations, we reveal the emergence of complex oscillatory behavior, including multi-year periodic cycles, transitions to heteroclinic cycle and chaotic oscillations. Our results highlight how periodic modulation can destabilize stable periodic trajectories, generate novel attractors, and give rise to coexistence mechanisms not present in autonomous systems. These findings contribute to the understanding of biodiversity maintenance under realistic, temporally varying ecological conditions.
   
\end{abstract}

\vspace{0.5em}

\section{Introduction: \large modeling of seasons in dynamical systems}
Many ecological systems adapt to seasonality\index{seasonality} 
by changing their behavioral patterns
in the course of the year~\cite{Vettorazzi}. 
One of the ways to model seasonality in dynamical systems 
is to introduce time-dependent system parameters;
in population dynamics these are e.g. the prey and birth rates.
In \cite{Taylor} it was demonstrated that applying seasonal modeling to ecological dynamical systems\index{ecological dynamical systems} can generate a rich variety 
of dynamics such as period-doubling bifurcation\index{period-doubling bifurcation}, multi-year cycles of different periods, parameter ranges with coexisting multi-year cycles\index{multi-year cycles} of the same or different periods, as well as 
quasi-periodicity and chaos\index{quasi-periodic}\index{chaos}.
Since seasonality is a roughly periodic phenomenon, below we assume the parameters of the system to be periodic functions of time.
For the case of long stationary seasons the latter can be approximated by
piecewise-constant step-functions\index{step-functions}, 
so that throughout each season the evolution is governed by 
a dynamical system with the constant set of parameters. 
Within this approximation, transitions between seasons are treated as 
instantaneous. Step-functions have been used in modeling the activity of the heart \cite{Atangana}; piecewise periodic functions are employed in the model of drug therapy of HIV~\cite{DeLeenheer}.

We restrict ourselves to the case when there are
only two  seasons, e.g. summer and winter. 
Let $P$ denote the overall period comprising both seasons.
Let $p_1$ and $p_2=P-p_1$ denote 
durations of, respectively, the first and the second seasons. 

Modeling seasonality by using piecewise constant functions means that some of the parameters become functions of time in the form 
\begin{equation}
\label{eq8}
\left\{\begin{array}{cc}
l_1,\quad & 0\leq t\,{\rm mod}\,P < p_1\, \\
l_2,\quad & p_1 \leq t\,{\rm mod}\,P\leq P.
\end{array}\right.
\end{equation}
{So far, the analysis of 
the equations of population dynamics} with parameters of the form (\ref{eq8}) has been carried out only in the context of “Seasonal Succession”: 
the species grow during the warmer months and die or form resting stages in winter. Mathematically, this means that part of the time dynamics is governed
e.g. by 
the Lotka-Volterra equations\index{Lotka-Volterra} whereas
during the rest of the time the amount of every species is decaying exponentially. 
In such systems, the main interest was in proving the stability of equilibrium points, periodic solutions, or heteroclinic cycles \cite{Xie2021a,Xie2021b,Xie2023}.

There is a significant difference between regular autonomous dynamical systems and dynamical systems with time-dependent parameters\index{time-dependent parameters} such as in (\ref{eq8}). First, the latter systems are non-smooth, hence not all results,
obtained for smooth dynamical systems, stay valid. 
In addition, the explicit dependence of the right-hand sides of the equations on time 
is equivalent to introducing an additional cyclic variable, 
which can make the dynamics more involved.

One can expect a significant change in the dynamics in such systems compared to their autonomous counterparts with constant parameters.

For example, below we shall see that in the considered exemplary 
autonomous systems  
chaotic behavior is impossible, but by making some of the parameters piecewise 
constant, we may, in principle, obtain chaos. 

In the literature, the phenomenon of switching between dynamical systems usually appears in the context of control theory \cite{Chen2024}. A typical example is a system of the form:
\begin{equation}
\label{eq9}
 \dot{x}=\left\{
\begin{aligned}
   f_1(x) \ , \ x\in D \\
   f_2(x) \ , \ x\notin D
\end{aligned}
\right.
\end{equation}
Within this setup, switching between dynamical systems does not directly depend on time but on the 
instantaneous state of the system: on the values of the system
variables. Systems of the type of (\ref{eq9}) are autonomous while systems of the type:
\begin{equation}
\label{eq10}
\dot{x}=\left\{
\begin{aligned}
    f_1(x,t) \ , \ t\in S \\
    f_2(x,t) \ , \ t\notin S
\end{aligned}
\right.
\end{equation}
are non-autonomous.

Below, in Section 2 we present the models to be studied and summarize their properties. Section 3 describes the extensions and details the results of our analysis. Section 4 provides a discussion of the findings, and Section 5 contains concluding remarks.
\section{Materials and methods: \large Multispecies competitive models}
\label{subsec:1.1}
\subsection{General Lotka-Volterra model}
\label{sub_gen_LV}
The general Lotka-Volterra model which describes temporal evolution
in a set of $n$ biological species competing for common resources, 
is a system of ordinary differential equations (ODEs)\index{ode}

\begin{equation}
\frac{dN_i(t)}{dt} =r_i N_i(t)\left[1-\sum_{j=1}^n \alpha_{ij}N_j (t)\right]
;\quad N_i\geq 0;\quad i=1,\ldots,n.
\label{gen_LV}
\end{equation}
Here $N_i(t)$ is the suitably normalized number of individuals of the $i$-th  species at time $t$, $r_i$ is the intrinsic growth rate of the $i$-th species, and $ \alpha_{ij}$ are the interaction coefficients that measure the extent to which the $j$-th species affects the growth rate of the $i$-th species \cite{CantrellCosner}.
For $n\geq 3$,  Eq.(\ref{gen_LV}) features rich dynamics, including multiple attractors, heteroclinic networks, and chaos~\cite{Postlethwaite2022,Pikovsky,Bayliss}. For $n=3$, the system (\ref{gen_LV}), after a proper rescaling, takes the form: 
\begin{eqnarray}
\label{eq2}
\dot{N_1}&=&r_1N_1\,(1-N_1-\alpha_1 N_2-\beta_1 N_3)\nonumber\\
\dot{N_2}&=&r_2N_2\,(1-N_2-\alpha_2 N_1-\beta_2 N_3)\\ 
\dot{N_3}&=&r_3N_3\,(1-N_3-\alpha_3 N_1- \beta_3 N_2)\nonumber
\end{eqnarray}
where the dot denotes the temporal derivative.
 
In \cite{Coste} it has been shown that under non-equal growth rates $r_i$ this system can have three kinds of attractors: a coexistence state of equilibrium (none of the variables equals zero), a limit cycle, and a heteroclinic cycle. 
The special case  $r_1=r_2=r_3=r$ is called the May-Leonard system~\cite{MayLeonard}\index{May-Leonard}. 
For this system, it was proven \cite{Chi} that for parameters that satisfy the relation 
$$(1-\beta_1)(1-\beta_2)(1-\beta_3)=(\alpha_1-1)(\alpha_2-1)(\alpha_3-1),\quad 0<\beta_i<1<\alpha_i,$$ 
the attractor can be a two-dimensional invariant manifold with a continuum of periodic solutions on it.

In the general case of (\ref{gen_LV}) with $n$ equal growth rates (called the \textit{general} May-Leonard  system), it has been proved that the system can be reduced to a competitive Lotka-Volterra system with $n-1$ species~\cite{Coste}. This implies that chaos can occur neither in (\ref{gen_LV}) for $n=3$ nor in the general May-Leonard system for $n=4$, since the latter can be reduced to a competitive Lotka-Volterra system with $n=3$.

\subsection{Symmetric May-Leonard model}

Here we discuss in more detail the paradigmatic {\em symmetric} May-Leonard model~\cite{MayLeonard} 
of three competing species, setting $r_1=r_2=r_3=1$, 
taking the circulant competition matrix 
$\begin{bmatrix}
1 &\alpha & \beta \\
\beta & 1 & \alpha \\
\alpha & \beta & 1
\end{bmatrix}$\\
with $0<\beta <1< \alpha$,
and re-denoting the variables 
$N_1=u$, $N_2=v$ and $N_3=w$. The model becomes:
\begin{eqnarray}
\label{eq3}
   \dot{u}&=&u\,(1-u-\alpha v- \beta w)\nonumber\\
   \dot{v}&=&v\,(1-v-\alpha w- \beta u)\\ 
   \dot{w}&=&w\,(1-w-\alpha u- \beta v).\nonumber
\end{eqnarray}
Dynamics\footnote{In the literature, the system (\ref{eq3}) can also be found in the form
\begin{equation*}
\begin{cases}
    \dot{a}\, =& a\left(1 - \rho - (\sigma + \zeta) b + \zeta c\right),\\
    \dot{b}\, =& b\left((1 - \rho - (\sigma + \zeta) c + \zeta a\right),\\
    \dot{c}\, =& c\left((1 - \rho - (\sigma + \zeta) a + \zeta b\right),
\end{cases}
 \end{equation*}
with $\rho=a+b+c$ and $\sigma,\zeta$ being the growth rates \cite{Postlethwaite2017}.} is considered in the invariant octant $u,v,w\geq 0$.
 The system (\ref{eq3}) has five equilibrium points: $(0,0,0)$, $(1,0,0)$, $(0,1,0)$, $(0,0,1)$, and the ``coexistence point"
\begin{equation}
    u=v=w=\frac{1}{1+\alpha+\beta}.
    \label{eq4}
\end{equation}
The coexistence point is stable, if $\alpha+\beta<2$, and is oscillatory unstable otherwise.
The planes $u=0,\,v=0$ and $w=0$ are invariant manifolds. On the plane $w=0$, there exists the heteroclinic trajectory leading from $(1,0,0)$ to $(0,1,0)$. Similarly, two more heteroclinic trajectories exist: from $(0,1,0)$ to $(0,0,1)$ on the plane $u=0$ and from $(0,0,1)$ to $(1,0,0)$ on the plane $v=0$. Thus, for any $\beta<1<\alpha$, there exists the robust heteroclinic cycle $(1,0,0)\to(0,1,0)\to(0,0,1)\to(1,0,0)$. Stability analysis\index{stability analysis} shows that the heteroclinic cycle\index{heteroclinic cycle} is attracting if $\alpha+\beta>2$, and is repelling if $\alpha+\beta<2$.
\noindent Therefore, the coexistence point is an attractor at $\alpha+\beta<2$, and the heteroclinic cycle is an attractor at $\alpha+\beta>2$. 

Dynamics at $\alpha+\beta=2$ is quite non-standard. 
At this combination of parameters, the plane $u+v+w=1$ is invariant and attracts all solutions that start outside this plane. Furthermore, dynamics 
is conservative: the combination 
\begin{equation}
A=\frac{u\,v\,w}{(u+v+w)^3}
\label{eq5}
\end{equation}
is a constant of motion which we define as the coexistence index.
Therefore, the attracting plane $u+v+w=1$ is foliated into closed phase trajectories. This continuum of periodic motions is parameterized by the
value of $A$ which, on that plane, equals $u\,v\,w$  and 
assumes the values from the interval $0\leq A\leq 1/27$.
 
The maximal value $A=1/27$ corresponds to the coexistence equilibrium $u=v=w=1/3$; the minimal value $A=0$ is reached at all other points of equilibrium and at the heteroclinic trajectories.

\begin{figure}
\centering
\includegraphics[width=0.5\linewidth]{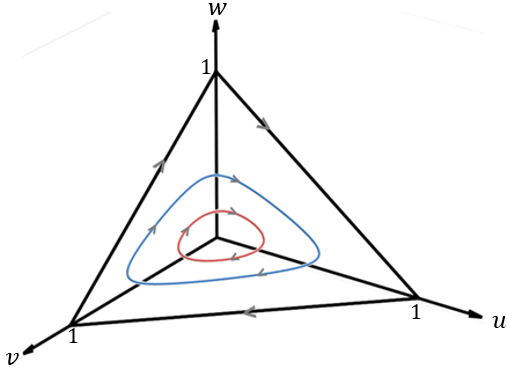}
\caption{\label{figure 1}\textbf{Phase portrait on the invariant simplex for $\alpha+\beta=2$.} On the plane $u+v+w=1$ (population fractions), the boundary triangle forms a heteroclinic cycle connecting the three corner equilibria. Trajectories in the interior are closed periodic orbits; two examples with different values of the conserved first integral $A$ are shown, illustrating that for $\alpha+\beta=2$ the dynamics on the simplex is conservative.}
\end{figure}

Systems with cyclic competition are often called Rock-Paper-Scissors systems due to their apparent similarity to the game: Scissors cut Paper, Paper wraps Rock, Rock blunts Scissors. In (\ref{eq3}), $v$ dominates $u$, $w$ dominates $v$, and $u$, in its turn, dominates $w$. In such systems, there is usually a heteroclinic cycle of the kind sketched in Fig.~\ref{figure 1}. 

Systems of this kind can be extended to games with more ``players" 
(species). For example, the game with 5 species called ``Rock–Paper–Scissors–Lizard–Spock" with the same cyclic rules as in the classical RPS game was studied in \cite{Postlethwaite2022}. In this variant of the game, the appearance of heteroclinic networks is possible, since there are species with ``non-direct interaction". 

\subsection{Specific aspects of dynamics}
\subsubsection*{Family of time-periodic solutions in the symmetric May-Leonard model}
\begin{figure}[ht]
    \centering
    \includegraphics[width=0.7\textwidth]{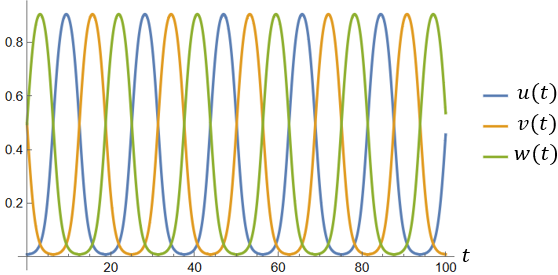}
    \caption{\label{figure2}\textbf{The species densities as functions of time: $u(t)$(blue), $v(t)$(orange) and $w(t)$(green) with the coexistence index $A=2\times10^{-3}$.}}
\end{figure}

\begin{figure}[htbp]
\centering
\includegraphics[width=0.5\linewidth]{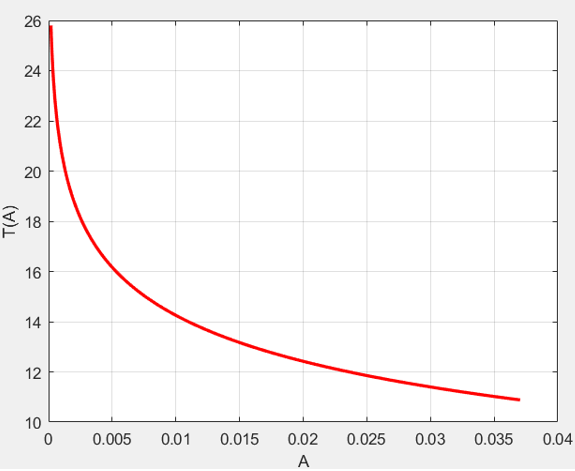}
\caption{\label{figure3} \textbf{Period $T$ versus \textcolor{black}{the coexistence index} $A$ (conservative case).}
The oscillation period tends to infinity as $A\to 0$ (near the heteroclinic cycle) and approaches to $2\pi \sqrt{3}$ as $A$ increases toward the coexistence point.}
\end{figure}

Consider Eq.(\ref{eq3}) at $\alpha=2$, $\beta=0$:
\begin{subequations}
\label{eq11}
  \begin{equation}  
\dot{u}=u\,(1-u-2v)
\label{eq11a}
\end{equation}
\begin{equation} 
\dot{v}=v\,(1-v-2w)
\label{eq11b}
\end{equation}
\begin{equation}
\dot{w}=w\,(1-w-2u)
\label{eq11c}
\end{equation}
\end{subequations}
Note that this May-Leonard model is symmetric with respect to cyclic permutation of the variables.
On the attracting manifold $u+v+w=1$, the closed trajectory 
$u\,v\,w=A$ is described by the equation \( uv(1-u-v)=A\), hence
\begin{equation}
v=\frac{u-u^2 \pm \sqrt{(u-u^2)^2-4Au}}{2u}.
\label{eq12}
\end{equation} 
On substituting expression (\ref{eq12}) 
into Eq.(\ref{eq11}), we find that
\begin{equation}
\dot{u}=\pm \sqrt{P(u)},\;P(u)=(u-u^2)^2-4Au,\; 0<u<1.
\label{eq13}
\end{equation}
In (\ref{eq13}), 
the positive and negative signs correspond, respectively, to the epochs of the growth and decrease of $u(t)$. The solution of (\ref{eq13}) with the initial condition $u(0)=u_1$, $0<u_1<1$ is
\begin{equation}
t(u) = \frac{2}{\sqrt{u_2 (u_3 - u_1)}} \, \text{F}\left( \arcsin \left( \sqrt{\frac{\displaystyle  \frac{1}{u} - \frac{1}{u_1}}{\displaystyle\frac{1}{u_2} - \frac{1}{u_1}} } \right), \frac{u_3 (u_2 - u_1)}{u_2 (u_3 - u_1)} \right).
\label{eq14}
\end{equation}
where $\text{F}(...)$ denotes the elliptic amplitude.
The inverse  function is
\begin{equation}
u(t)= \frac{1}{\displaystyle\frac{1}{u_1}+(\frac{1}{u_2}-\frac{1}{u_1})\cdot {\rm sn}^2\left[\frac{t}{2}\sqrt{u_2(u_3-u_1)},\frac{u_3(u_2-u_1)}{u_2(u_3-u_1)}\right]}
\label{eq15}
\end{equation}
 where $u_1 < u_2 <u_3$ are the positive roots of the polynomial $P(u)$ (the fourth root is zero), and
sn denotes the elliptic sine function.
Analyzing the properties of the polynomial $P(u)$ with $0<A<1/27$, we conclude that \(u_1 , u_2<1 , u_3>1\),
so that $u(t)$  varies between $u_1$ and $u_2$, because the expression under the square root should be positive and the minimal and maximal values of $u$  correspond to $\dot{u}=0$. Due to
the permutational symmetry of the problem, the solutions $v(t)$, $w(t)$ are the replicas of $u(t)$ shifted by $T/3$:
\begin{equation}
v(t)=u\left(t-\frac{T}{3}\right) \ , \  w(t)=u\left(t+\frac{T}{3}\right),
\label{eq16}
\end{equation}
 where  \(T\) is the period of (\ref{eq16}), which can be found as
$T=2t(u_2)$, i.e.,
\begin{equation}
T= \frac{4}{u_2(u_3-u_1)}K\left(\frac{u_3(u_2-u_1)}{u_2(u_3-u_1)}\right),
\label{eq17}
\end{equation}
where $K$ is the complete elliptic integral of the first kind
(see Figs~\ref{figure2},\ref{figure3}).
 The limits of period at the boundaries of  the interval $(0, 1/27)$ are: 
\begin{equation}
\label{eq18}
\lim_{A\to 0} T = 2\pi K(1)=\infty, \;
\lim_{A\to 1/27} T = 2\pi \sqrt{3}.
\end{equation}

\subsection{A simple dissipative model}
Consider now the Lotka-Volterra system
\begin{eqnarray}
\label{eq19}
u_t &= u(\gamma - u - \alpha v),\nonumber\\  
v_t &= v(1 - v - \alpha w)\\  
w_t &= w(1 - w - \alpha u)\nonumber  
\end{eqnarray}
with the growth rate of the first species different from that of two other species.   

\subsubsection{Linear instability of the coexistence state.}

The coexistence solution is
\begin{equation}
\label{eq20}
u_*=\frac{\gamma -\alpha +\alpha ^2}{1+\alpha ^3},\; 
v_*=\frac{1 -\alpha +\alpha ^2 \gamma}{1+\alpha ^3},\; 
w_*=\frac{1 -\alpha \gamma +\alpha ^2}{1+\alpha ^3}.
\end{equation}

We assume that
\begin{equation}
\label{eq21} 
\alpha >1,\; \frac{\alpha -1}{\alpha ^2}<\gamma <\frac{\alpha ^2 +1}{\alpha},
\end{equation}
hence $u_*$, $v_*$, $w_*$ are positive. 

Stability of the coexistence state is determined by the linearization of (\ref{eq19}):

three eigenvalues $\sigma$ of the Jacobian matrix solve the equation
\begin{equation}
\label{eq23}
\sigma ^3+(u_*+v_*+w_*)\sigma ^2+(u_*v_*+v_*w_*+w_*u_*)\sigma+(1+\alpha ^3)u_*v_*w_*=0
\end{equation}
At the stability boundary the Hopf bifurcation\index{Hopf bifurcation} takes place: two complex conjugate eigenvalues are purely imaginary, $\sigma=\pm i\omega$ for a real $\omega$.
This, after the substitution of the coordinate values, 

determines the bifurcation curve on the parameter plane $(\gamma,\alpha)$:
\begin{equation}
\label{eq25}\begin{split}
\alpha (1 - \alpha - \alpha^3) \gamma^3 
+ (-2 + 3\alpha - 5\alpha^2 + 6\alpha^3 + \alpha^5 - \alpha^6) \gamma^2\\ 
+ (-4 + 7\alpha - 11\alpha^2 + 5\alpha^3 - 7\alpha^4 + \alpha^5 + \alpha^7) 
\gamma\\ 
+ (1 - \alpha)^2 (-2 + \alpha - 3\alpha^2 - \alpha^3 - \alpha^4) = 0
\end{split}
\end{equation}
The curve of Hopf bifurcation splits the parameter plane into the region 
where the coexistence point is the attractor, and the region where it is unstable.

\subsubsection{Stability of the heteroclinic cycle}

System (\ref{eq19}) has three saddle points $(\gamma,0,0)$, $(0,1,0)$, $(0,0,1)$, which are connected by heteroclinic trajectories forming a heteroclinic cycle.

Consider a trajectory that is very close to the heteroclinic cycle.
The time needed 
to complete a ``cycle" (start close to one of the saddle points and then come back to its neighborhood) can be divided into 6 intervals: 3 intervals where the solution is close to one of the saddles, with means that only one of the variables is $O(1)$, and 3 intervals on each of the heteroclinic trajectories, where two of the variables are $O(1)$. We will see that the major contribution comes from the time spent near the saddles.

Consider the trajectory starting between the vicinities of saddle points $(0,0,1)$ and $(\gamma,0,0)$, where $u$ and $w$ are $O(1)$, while $v=\epsilon \ll 1$.
Near the saddle point $(\gamma,0,0)$, the variables are governed by the linear equations:
$$\tilde{u}_t=-\gamma \tilde{u}-\alpha \gamma v, v_t=v, w_t=(1-\alpha \gamma )w.$$
Here $\tilde{u}=u-\gamma$.
The time $t_1$ spent near that point depends on $\epsilon$: during that time, the variable $v$ grows from $\epsilon$ to $O(1)$. Since $\Tilde{v} \sim \epsilon e^t$, we have $\epsilon e^{t_1} \sim 1$, so $t_1 \sim -\ln\epsilon$. 
During that time, $w$ decreases from $O(1)$ to $\Tilde{w_1} \sim e^{(1-\alpha \gamma)t_1}$. 

When $v$ becomes $O(1)$, the linear equation for $v$ has to be corrected.
Between the saddles $(\gamma,0,0)$ and $(0,1,0)$, where both $v$ and $1-v$ are $O(1)$, the governing equation for $v$ is $v_t=v(1-v)$ and the solution for $v$ corresponds to the heteroclinic trajectory in the plane $w=0$,
\begin{equation}
\label{eq26}
v(t)=\frac{1}{1+C_1e^{-t}},
\end{equation}
where $C_1$ is a constant depending on the initial conditions. Because for $v=O(1)$ and $1-v=O(1)$, the velocity $v_t=O(1)$, the duration of the motion on that fragment of trajectory is $O(1)$. Hence, it is much smaller than $t_1$, and the variable $w$ does not change its order during that interval, hence the change of the variable $w$ during that time can be neglected.
Near the point $(0,1,0)$, the evolution of variables is governed by the system of equations $$u_t=(\gamma-\alpha)u, \ \tilde{v_t}=-\tilde{v}-\alpha w,\; w_t=w,$$
where $\tilde{v}=v-1$.
During the time $t_2$ spent near that point, $w$ growth from $\Tilde{w_1}$ up to $O(1)$, hence, similarly to the calculation of $t_1$ we get: $t_2\sim -\ln(\Tilde{w_1})=-(\alpha \gamma -1)\ln\epsilon$.
During that time $t_2$, $u$ decreases from $O(1)$ to $u_2\sim e^{(\gamma-\alpha)t_2}=\epsilon^{(\alpha-\gamma)(\alpha \gamma -1)}$.

Between the points $(0,1,0)$ and $(0,0,1)$ the governing equation for $w$ is $w_t=w(1-w)$,
and the heteroclinic solution on the plane $u=0$  is
\begin{equation}
\label{eq27}
w(t)=\frac{1}{1+C_2e^{-t}},
\end{equation}
where $C_2$ is a constant depending on the initial conditions. Again, the duration of the time interval where $w=O(1)$ and $1-w=O(1)$ is $O(1)$, and the variable $u$ does not change its order during that interval.

Near the point $(0,0,1)$, the evolution of variables is governed by the system of equations,
$$u_t=\gamma u, \ v_t=(1-\alpha)v, \ \tilde{w_t}=-\tilde{w}-\alpha \tilde{u},$$
where $\tilde{w}=w-1$.
The variable $u$ grows from 
$\Tilde{u_2}$ till $O(1)$ during the time $t_3$: $t_3\sim -\frac{1}{\gamma}\ln(\Tilde{u_2})=-\frac{(\alpha -\gamma)(\alpha \gamma -1)}{\gamma}\ln\epsilon$. \\
During that time $v$ decreases from $O(1)$ to $\Tilde{v_3}\sim e^{-(\alpha-1)t_3}=\epsilon^{(\alpha-1)(\alpha-\gamma)(\alpha \gamma-1)/\gamma}$. \\
Finally, between the points $(0,0,1)$ and $(\gamma,0,0)$ the governing equation for $u$ is $u_t=u(\gamma-u)$ and the solution for $u$ corresponds to the heteroclinic trajectory in the plane $v=0$,
\begin{equation}
\label{eq28}
u(t)=\frac{\gamma}{1-e^{\gamma(C_3-t)}}
\end{equation}
where $C_3$ is a constant depending on the initial conditions.
The variable $v$ does not change its order on that interval.

\begin{figure}[htbp]
\centering
\includegraphics[width=1.0\linewidth]{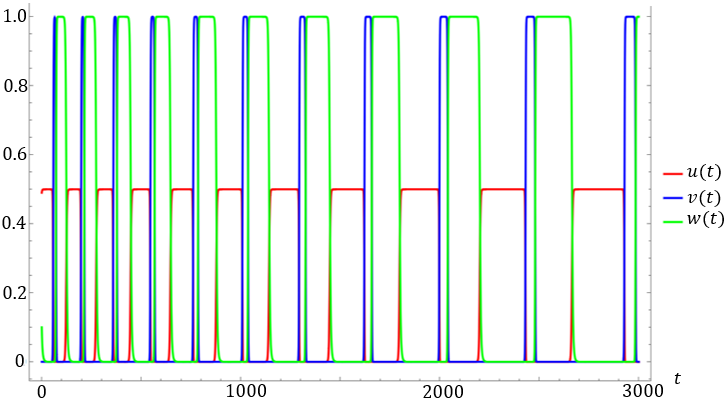}
\caption{\label{figure4} \textbf{Approach to the heteroclinic cycle.}
An orbit on the invariant simplex spends progressively longer near the saddle equilibria and travels along the heteroclinic connections, illustrating the slow-down toward the heteroclinic cycle.}
\end{figure}

We find that during the motion in the vicinity of the heteroclinic cycle $(\gamma,0,0)\rightarrow (0,1,0)\rightarrow (0,0,1)\rightarrow (\gamma,0,0)$ the order of the variable $v$ is changed from $\epsilon$ to $\epsilon^c$ where 
$$c=(\alpha-1)(\alpha-\gamma)(\alpha \gamma -1)/\gamma.$$ 
When $c>1$, the heteroclinic cycle is an attractor since a solution that starts close to the heteroclinic cycle becomes even closer to it. When $c<1$, heteroclinic cycle is a repeller. Note that if we start from a different point, we obtain the same coefficient $c$ (with different order of multiplications). The critical case $c=1$ divides the plane $(\gamma,\alpha)$ by the curve:
\begin{equation}
\label{eq29}
(\alpha-1)(\alpha \gamma -1)(\alpha -\gamma)-\gamma=0
\end{equation}
The meaning of the latter expression is as follows. In accordance to the general theory, stability of the heteroclinic cycle requires that compression near the cycle prevails over expansion: product of all positive Jacobian eigenvalues for the participating saddles is smaller than the absolute value of the product of all closest to zero negative Jacobian eigenvalues. Equality between these two products, in terms of the coefficients of the equations, results in the Eq.(\ref{eq29}).
Near the heteroclinic cycle, for each solution long epochs of nearly quiescent hovering over one of the saddle points alternate with relatively fast passages between the saddle points. The time of the latter passage is negligible compared to the time of slow creeping near the saddle points. Since $c>1$, during the next "cycle" the time $t_1$ spent near the point $(\gamma,0,0)$ will be $-c\ln\epsilon$, so we conclude that the time that the solution spends near each saddle increases
with each new turn around the cycle. This behavior can be seen in Fig.~\ref{figure4}.
\subsubsection{Attractors in the system of ODEs}

\begin{figure}[htbp]
\centering
\includegraphics[width=0.7\linewidth]{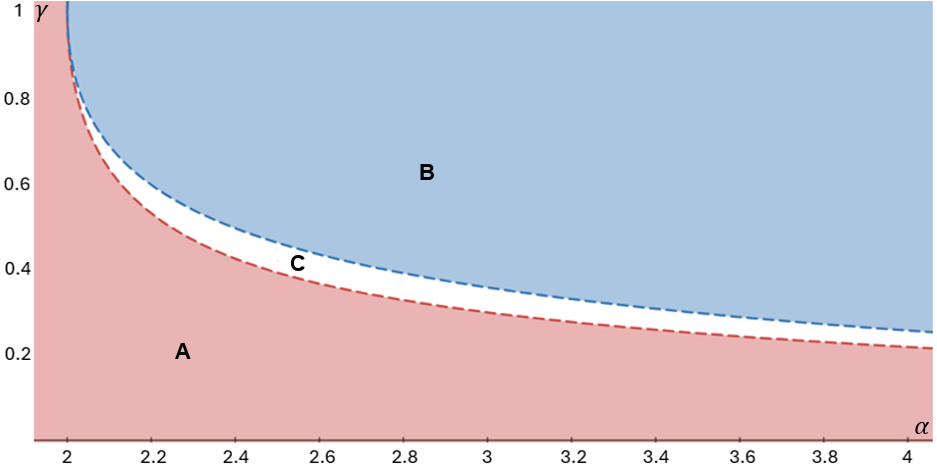}
\caption{\label{figure5}Parameters (interaction coefficient and growth rate) plane $(\alpha,\gamma)$: attractor regimes and bifurcation curves.
Colored shading distinguishes three regions: 
\textbf{A}-stable coexistence equilibrium; 
\textbf{B}-stable heteroclinic cycle; 
\textbf{C}-stable limit cycle (between \textbf{A} and \textbf{B}). 
The boundary of \textbf{C} is formed by the Hopf bifurcation curve (where the limit cycle is born) and the curve on which the heteroclinic cycle becomes stable; these curves are tangent at $(\alpha,\gamma)=(2,1)$ (the symmetric May--Leonard case).} 

\end{figure}

Equations (\ref{eq25}),(\ref{eq29}) together with inequalities (\ref{eq21}) split the parameter plane $(\alpha,\gamma)$ into three regions: the region $A$ where the attractor is the equilibrium, the region $B$ where the attractor is the heteroclinic cycle and the region $C$ between them. 
In \cite{Coste} it was shown that for the systems in the form of (\ref{eq19}) only three possible attractors were possible: a state of equilibrium, a heteroclinic cycle and a limit cycle\index{limit cycle}. In accordance with this statement,
the attractor in the region $C$ should be the limit cycle.
Indeed,  on the state diagram, shown in Fig.~\ref{figure5}, the borders of the region $C$ are the line of the Hopf bifurcation (in which the limit cycle is born) and the line  (\ref{eq25}) on which the heteroclinic cycle acquires stability.
Two bifurcation curves on the parameter plane are mutually tangent in the point $(\alpha =2,\gamma=1)$ which is the symmetric May-Leonard model (\ref{eq11}).

\begin{figure}[htbp]
\centering
\includegraphics[width=0.7\linewidth]{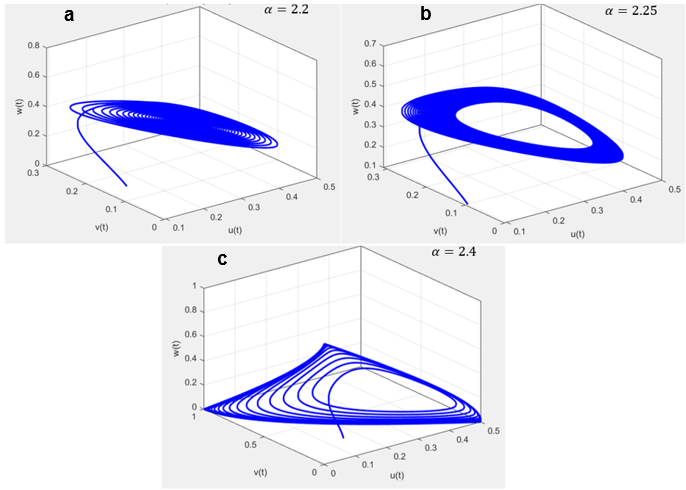}
\caption{\label{figure6} \textbf{Phase trajectories at fixed growth rate $\gamma$ as the interaction coefficient $\alpha$ varies.}
With the same initial condition, the flow shows (a) convergence to the coexistence equilibrium for $\alpha$ below the Hopf bifurcation, (b) a stable limit cycle just beyond Hopf, and (c) an attracting heteroclinic cycle once the heteroclinic threshold is crossed.
}
\end{figure}

From now on we fix $\gamma=0.5$ and treat $\alpha$ as the control parameter of the system.  For $\alpha=2.24:=\alpha_{Hopf}$ the supercritical Hopf bifurcation occurs and the stable limit cycle appears, while for $\alpha=2.38:=\alpha_{het}$ the heteroclinic bifurcation takes place and the heteroclinic cycle acquires stability. 
These changes in the phase space\index{phase space} are presented in Fig.~\ref{figure6}.

\newpage
\section{Results}
\subsection*{Alternating May-Leonard and Lotka-Volterra dynamics}
This section is motivated by modeling seasonal variation in ecological systems.
\subsection{From the set of differential equations to the stroboscopic map}
Consider the system
\begin{equation}
\label{eq30}
\dot x(t) =
\begin{cases}
f_{1}(x), & 0 \le t\; \mathrm{mod}(p_{1}+p_{2})\, \le p_{1},\\
f_{2}(x), & p_{1} < t\; \mathrm{mod}(p_{1}+p_{2}) < p_{1}+p_{2},
\end{cases}
\end{equation}

with $f_1,f_2$ being functions of the type (\ref{eq8}), specifically, of
the Lotka-Volterra type or the May-Leonard type.
We choose $x$ to be the 3-species vector $(u,v,w)^T$ that obeys the system 
\begin{equation}\label{eq31}
\begin{alignedat}{2}
u_t \;&&={}&\;u\bigl(\gamma(t)-u-\alpha(t)\,v\bigr),\\
v_t \;&&={}&\;v\bigl(1-v-\alpha(t)\,w\bigr),\\
w_t \;&&={}&\;w\bigl(1-w-\alpha(t)\,u\bigr).
\end{alignedat}
\end{equation}
where $\alpha,\gamma$ are piecewise constant functions of time:
\begin{equation}\label{eq32}
(\alpha(t),\gamma(t)) =
\begin{cases}
(\alpha_{1},\gamma_{1}), & 0 \le t\; \mathrm{mod}(p_{1}+p_{2}) \le p_{1},\\
(\alpha_{2},\gamma_{2}), & p_{1} < t\: \mathrm{mod}(p_{1}+p_{2}) < p_{1}+p_{2},
\end{cases}
\end{equation}

The system (\ref{eq31}) together with the rule (\ref{eq32}) has the form (\ref{eq30}). As time goes on, the parameter values are ``jumping'' periodically between the set $(\alpha_1,\gamma_1)$ and the set $(\alpha_2,\gamma_2)$. To study dynamics, we choose 
the combination of $(\alpha_1,\gamma_1)=(2,1)$ which yields the symmetric May-Leonard model (\ref{eq11}) and $(\alpha_2,\gamma_2)=(2.26,0.5)$ which renders the Lotka-Volterra model (\ref{eq19}) with the limit cycle of period 16.998849 as the global attractor.
Accordingly, we denote the time interval $0\leq t \leq p_1 $ as conservative season and the interval $ p_1< t < p_1+p_2 $ as dissipative season.

\begin{figure}[htbp]
\centering
\includegraphics[width=0.95\linewidth]{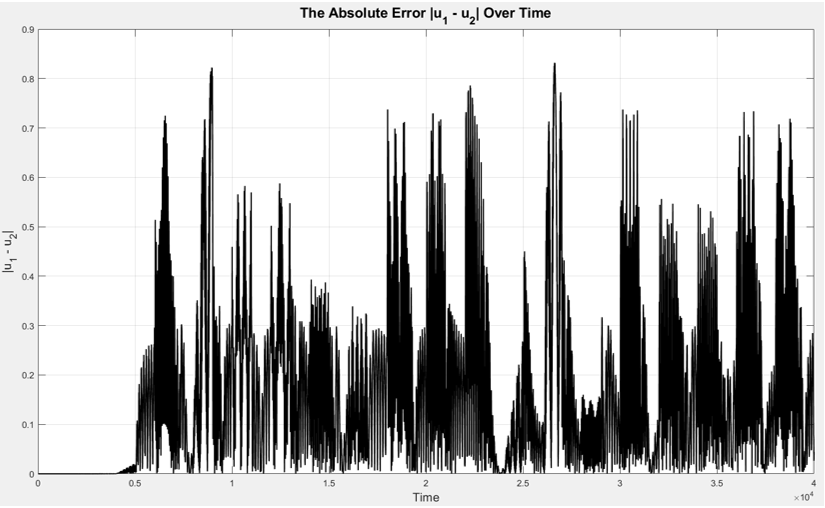}
\caption{\label{figure7}\textbf{Sensitivity to initial conditions under seasonal switching.}
Two solutions of the species density $u(t)$ starting $10^{-9}$ apart separate markedly after several switching cycles with the seasons length $p_1=p_2=1000$ fixed, demonstrating strong sensitivity to initial conditions.}
\end{figure}

Introduction of periodic parameter switching results in novel dynamical effects, impossible in a setup with fixed parameters. One of these effects reminds sensitive dependence on initial conditions, known as a hallmark of chaotic dynamics~\cite{Ott_book}.
In Fig.~\ref{figure7}, temporal evolution of the difference between two initially rather close solutions (difference of $10^{-9}$ in $u$ and $v$) is plotted. It can be seen that at the beginning the solutions start very close to each other, but after enough time, the difference between them reaches its maximal value $\sim 0.85$.

\begin{figure}[htbp]
\centering
\includegraphics[width=1\linewidth]{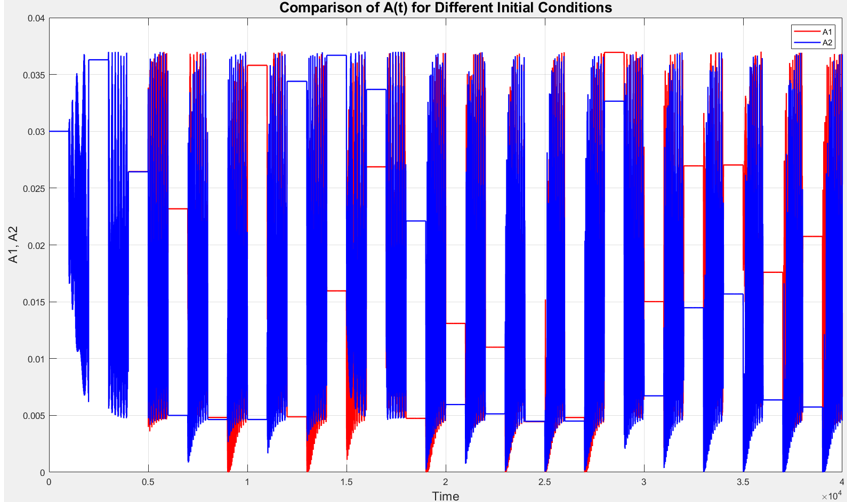}
\caption{\label{figure8}\textbf{Energies at seasonal boundaries for two nearby solutions.}
Within each conservative season the stroboscopic coexistence index $A$ is constant (plateaus). After about three full forcing cycles $(p_1+p_2)$, the two solutions enter conservative seasons with different $A$ (coexistence index)} values, so their plateaus no longer coincide.
\end{figure}

Since the solution of the conservative system\index{conservative systems} is defined by its coexistence index $A$, it is instructive to follow the evolution of the coexistence index values. In Fig.~\ref{figure8} the coexistence index levels have separated across conservative seasons\index{conservative systems}, while within any given conservative season $A$ remains constant and changes only during the dissipative stage.

In order to understand the origin of this effect, we introduce a new tool:
sequence of states 
at certain specially chosen values of $t$.
If $p_2$ is sufficiently big, the system has enough time to approach the limit cycle from any generic initial state, hence we expect that after the completion of every period $p_1+p_2$ 
the system is, in a sense, close to that cycle.
In the symmetric May-Leonard case all trajectories are 2-d manifold ${u\,v\,w}/{(u+v+w)^3}=A$, and the location on a particular closed trajectory depends on the initial conditions in the phase space: the initial value of the  coexistence index $A$. Note that during the dissipative season the value of $A$ is not conserved and oscillates along the limit cycle; at the employed values of $\alpha$ and $\gamma$ the oscillation range is $0.00464<A<0.036988$. The period of the latter cycle is the ``internal'' characteristics of the Lotka-Volterra dynamics; in general, it is incommensurate with ``externally'' imposed season durations $p_1$ and $p_2$. 
Therefore, after every full period $p_1+p_2$ the system starts its next conservative stage with a next value of $A$, and is, as a result, attracted to a different periodic solution on the two-dimensional manifold. By following in the phase space the sequence of endpoints of conservative seasons (or, alternatively, the sequence of endpoints of dissipative seasons), we obtain a mapping of the three-dimensional phase space $(u,v,w)$ into itself. It is important to note that this mapping is smooth although the original system of differential equations is non-smooth: it is a product of two smooth mappings which interrelate the starting point and the endpoint of the season (one mapping for the dissipative season, another mapping for the conservative one). 

We observe the system at the time instants  
$\{t_n\}=\{(p_1+p_2)n\},n=1,2,\ldots$ (Here and henceforth, $n$ denotes the season/time index; the number of species is fixed (three) and is not denoted by $n$)
and obtain in this way a stroboscopic map
$x_{n+1}(t_{n+1}) = f(x_n(t_n))$,
where $x_n$ is the triple of coordinates $\{u,v,w\}$ at $t_n$.
 
To understand the origin of complicated temporal dynamics, it makes sense to move on the parameter plane spanned by the season durations $p_1$ and $p_2$ along the lines parallel to the axes. Below, we view the behavior of the system when one of the durations: $p_1$  of the conservative stage or $p_2$ of the dissipative one stays frozen, and treat the remaining parameter as the free (and, in a sense, small) one.

\subsection{Fixed duration of the dissipative stage}
\begin{figure}[htbp]
\centering
\includegraphics[width=0.97\linewidth]{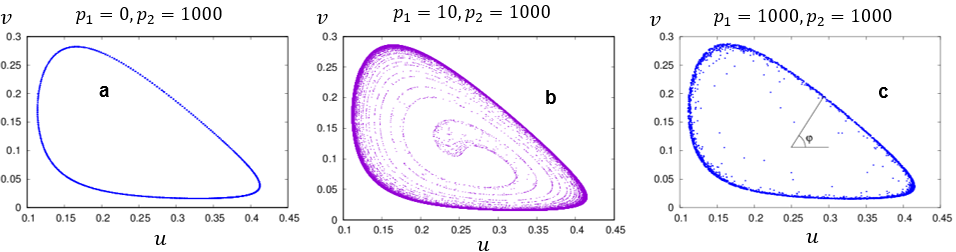}
\caption{\label{strobo_3} Two-dimensional projections ($u$ vs $v$) of the attractors
of the stroboscopic mapping at fixed duration $p_2=1000$ of the dissipative season, combined with vanishing (a), moderate (b)
and balanced (c) duration $p_1$ of the conservative season.}
\end{figure}

To start with, we fix the duration $p_2$ at some value that distinctly exceeds the period of the dissipative limit cycle (we use for illustrations $p_2=1000$), 
and compare the return mappings at different values of $p_1$.
Typical examples of attracting sets are shown in Fig.~\ref{strobo_3}.
At $p_1=0$ (that is, under a vanishingly short conservative season)  
all ``probing points'' lie precisely on the limit cycle of the dissipative system\index{dissipative systems}.
At $p_1=10$ (a non-zero, remarkably short conservative stage)  
the points fill the layered folded curve that lies inside the cycle. 
Finally, at $p_1=1000$, when both seasons have equal durations,
the majority of points is located in the narrow stripe around the limit cycle; however, a closer look
shows that the attracting set is ``fuzzy'' and disordered: there are
also seldom irregular excursions inwards, far away from the curve of the limit cycle, that look like noisy outliers.

The full stroboscopic map\index{stroboscopic map} operates in the three-dimensional space, but 
in the asymptotics of short conservative seasons it is well approximated by the one-dimensional circle map\index{circle map}. In this situation, most of the time the motion is one-dimensional and occurs along the limit cycle; short conservative seasons act as a regular sequence of ``kicks'' that temporarily drive the system in the phase space away from the curve of the limit cycle. If the conservative stages are sufficiently short, the kicks can be viewed as weak perturbations of the periodic motion (the fact that these perturbations are conservative does not matter much). A motion along the limit cycle is adequately described in terms of the phase: the angular variable that parameterizes the closed curve in the phase space~\cite{Glass}. 

Operationally, near the limit cycle the phase can be recovered by usage of the kind of ``polar coordinates'': for a sufficiently large collection of the mapping points, 
we find the coordinates of the center of mass (e.g. $\langle u\rangle $), 
and for (any two out of three) coordinates calculate the angle $\varphi_n$ between the line that connects the center of mass with the point  $(u_i,v_i)$ and one of the coordinate axes, e.g.
$$\tan(\varphi_i)=\frac{v_i-\langle v\rangle }{u_i-\langle u\rangle}.$$
Reconstruction of the angular variable $\varphi$ is sketched in Fig.~\ref{strobo_3}(c).
By varying $i$, we observe the sequence of the iterations of $\varphi_i$: the desired circle mapping 
\begin{equation}
\label{eq33}
\varphi_{i+1}=g(\varphi_i)\  {\rm mod}(2\pi)
\end{equation}

In the fuzzy situations of ~\ref{strobo_3}(b,c) the value of the angle alone does not determine the evolution, and the circle map is not one-to-one; nevertheless, it helps to disclose certain dynamical features.

\begin{figure}[htbp]
\centering
\includegraphics[width=0.87\linewidth]{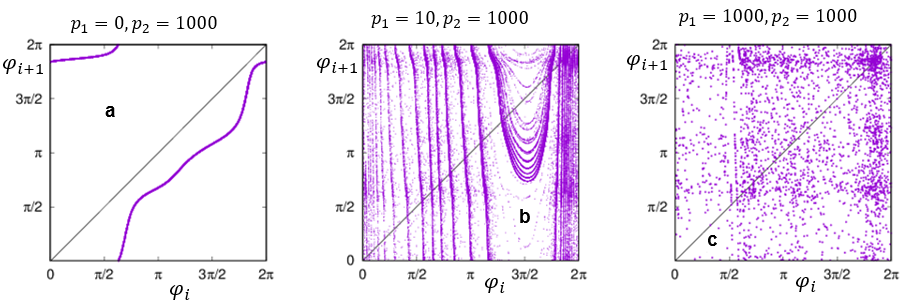}
\caption{\label{maps_3} Circle maps of the angle $\varphi_i$ between the line that connects the center of mass with the point  $(u_i,v_i)$ and one of the coordinate axes,  constructed from 
the stroboscopic mapping at fixed duration $p_2$ of the dissipative season,
and varying duration $p_1$ of the conservative season.
The values of $p_1$ in the panels (a),(b),(c) coincide with respective values 
from the panels of Fig.~\ref{strobo_3}.}
\end{figure}

Circle maps, recovered from the stroboscopic maps of
Fig.~\ref{strobo_3}, are presented in Fig.~\ref{maps_3}.

The plot of $g(\varphi)$ at $p_1=0$ in Fig.~\ref{maps_3}(a) features
the shape, characteristic for a circle map with irrational rotation number. The curve of the map is piecewise monotone
and does not intersect the ``bisecting line'' 
$\varphi_{i+1}=\varphi_i$. In this case there is no synchronization between two timescales of the system: the ``internal'' period of the limit cycle and the ``externally imposed'' full duration $p_1+p_2$.

At $p_1=10$ the plot of the circle map, shown in Fig.~\ref{maps_3}(b), is fuzzy, but still possesses a certain
pattern. There are a few dozens of ``branches'', each of them covering the whole vertical span of the map. Each of these branches intersects the bisecting line. Every intersection
delivers a fixed point $\varphi_{i+1}=\varphi_i$ of the circle map -- a periodic oscillation in the underlying system of ODEs. Since the absolute value of the slope of $g(\varphi)$ in these fixed points exceeds 1, all of them are unstable. Presence of the countable set of unstable periodic orbits is an indicator of chaotic dynamics~\cite{Ott_book}. Another notable feature of the map in Fig.~\ref{maps_3}(b) is the multi-branched layered local minimum in the right part. Such locally unimodal mappings are akin to the logistic mapping and provide the framework for the sequence of period-doubling bifurcations.

Finally, the plot of the circle map at $p_1=1000$, shown in Fig.~\ref{maps_3}(c),
bears no particular pattern. A multitude of points in the vicinity of 
the line $\varphi_{i+1}=\varphi_i$ indicates strong chaoticity. 

Apparently, the onset of chaos occurs during the transformation between
the left and the central panel of Fig.~\ref{maps_3}. Here, we briefly describe
the sequence of events and attracting states, encountered in the course of increase of $p_1$. A convenient characteristics is the observable  coexistence index $A_i={u\,v\,w}/{(u+v+w)^3}$ (or any other explicit function of the phase space coordinates) measured at the endpoint of the $i$-th dissipative season.
\begin{itemize}
\item In the range $0<p_1<p_{\rm syn}=0.091\ldots$ the shape of the map qualitatively reminds Fig.~\ref{maps_3}a. There is no synchronization between
the timescales of the problem. As $p_1$ is increased, the curve $g(\varphi)$
comes closer to the line $\varphi_{i+1}=\varphi_i$. Observable $A_i$ attains
a range of values.
\item At $p_1=p_{\rm syn}$ the curve $g(\varphi)$  becomes
tangent to the line $\varphi_{i+1}=\varphi_i$. This is the so-called
\textit{tangent} bifurcation: the structurally unstable neutral fixed point of
the circle map appears, and, under the further increase of $p_1$, decomposes
into \textit{two} fixed points: the stable and the unstable one. The stable
fixed point corresponds to the stable periodic solution in the system of ODEs.
The tangent bifurcation marks the onset of synchronization: in the parameter
space, the system enters the synchronization region -- the so-called Arnold 
tongue~\cite{Aquillen}\index{Arnold 
tongue}. Observable $A_i$ converges to the single constant value.
\item In the range $p_{\rm syn}<p_1<p_{\rm pd}=0.194\ldots$ the
stable fixed point of the circle map persists; oscillations in the system stay synchronized 
with the total duration of the seasons. Close to $p_{\rm pd}$, the mapping curve
becomes non-monotonic.
\item At $p_{\rm pd}$ the period-doubling bifurcation takes place in the map:
the fixed point loses stability, and a stable periodic orbit composed of two points:
$\varphi_{i+1}=g(\varphi_i),\;\varphi_{i+2}=g(\varphi_{i+1})=\varphi_{i}$
is born from it. Accordingly, the observable $A_i$ jumps between two alternating values. 
\item In the range $r_{\rm pd}<p_1<p_{\rm pd2}=0.23\ldots$ the mapping orbit
with period 2 is stable. At $p_{\rm pd2}$ the second period-doubling bifurcation occurs: the orbit with period 2 loses stability, and the stable orbit with period 4 is born. Just beyond $p_{\rm pd2}$, the observable $A_i$ jumps between \textit{four} alternating values. 
\item At $p_{\rm pd3}=0.236$ the further period-doubling bifurcation gives birth to the stable orbit with period~8.
\item At $p_{\rm pd4}=0.2377$ the next period-doubling bifurcation creates the stable orbit with period 16.
\item $\ldots$ At $p_{{\rm pd}j}$ the period-doubling bifurcation destabilizes the map orbit with period $2^{j-1}$ and creates the stable orbit with period $2^j$. The coexistence index $A$ jumps between \textit{$2^j$} alternating values.
\item At $p_1\geq 0.238$ we observe weak chaos in the dynamics of $A$. The plot of the circle map consists of several filled segments, some of them folded.
\item At $p_1\approx 0.27$ the local structure of the
circle map changes; now it has a shape of a narrow $M$-shaped band with a well-pronounced minimum, and consists of up to 3 thin layers.
\item At $p_1\approx 0.45$ the circle map is a narrow $M$-shaped band built of 5 thin layers.
\end{itemize}

\begin{figure}[htbp]
\centering
\includegraphics[width=0.8\linewidth]{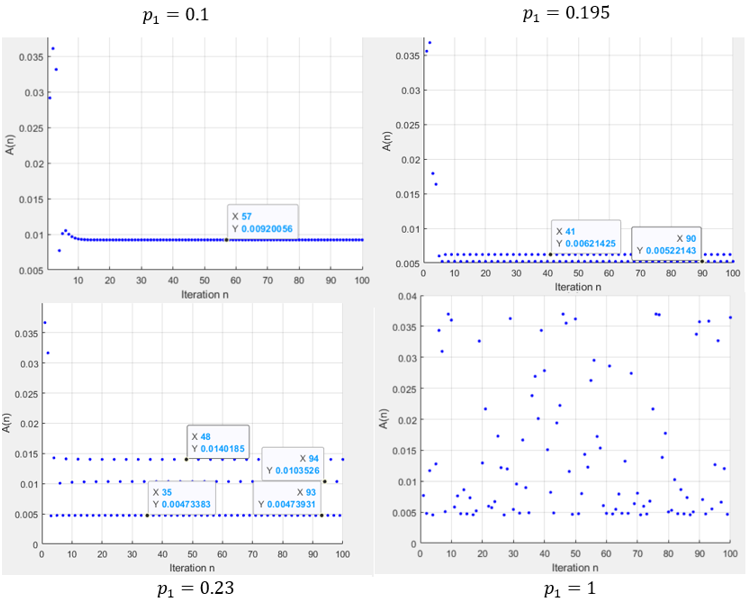}
\caption{\label{dissipative}
\textbf{Evolution of the coexistence index} $A_n$ for various conservative season durations $p_1$ at fixed dissipative season duration $p_2=1000$.}
Panels illustrate locking plateaus with one, two, or four levels as the $p_1$ increases, followed by irregular fluctuations indicating loss of phase locking.
\end{figure}

Altogether, the evolution of the circle map follows the celebrated scenario of period-doubling bifurcations: a sequence of intervals of $p_1$ along which the period length is progressively doubled. 
Each parameter interval of $p_1$ between two successive transitions is smaller than the preceding one by the factor that should asymptotically tend to the universal value 
$\delta_F=4.669201...$: the Feigenbaum constant~\cite{Feigenbaum}\index{Feigenbaum constant}. 
This behavior is common for generic one-parameter families of unimodal mappings, the most famous representative being the logistic map~\cite{May}. Like in all those families, beyond the first
threshold of chaos there is a set of windows of periodic behavior;
for example, close to $p_1=0.25$ the mapping again has a stable orbit of period 4 (which is different from the mentioned orbit born at $p_{\rm pd2}$ and destabilized at $p_{\rm pd3}$), but basically chaos prevails in the parameter space. For some of the numerical results See Fig. \ref{dissipative}.

\subsection{Fixed duration of conservative stage}
We consider the system (\ref{eq31}) with piecewise-constant coefficients in the right hand side (\ref{eq32}) under the assumption that the conservative season is sufficiently long, i.e., $p_1 \gg 1$. In this regime, the system has enough time to approach a periodic orbit characterized by a nearly conserved total population, so that $u + v + w \approx 1$ at the end of the conservative phase. \\
We first recall that for $p_2=0$, solutions are attracted to periodic orbits determined by the initial value of $A$. Here we study the effect of a short dissipative stage $p_2 \ll 1$. \\
The system is governed by equations
$\dot{u}=u(\gamma-u-\alpha v),\,\dot{v}=v(1-v-\alpha w),\,\dot{w}=w(1-w-\alpha u)$; during the dissipative stage the parameter values are $\alpha=2.26,\gamma=0.5$. Let the endpoint of the conservative season in the phase space be the state $(u_c,v_c,w_c)$. Since $p_2$ is short, expressions for the increments of the variables during the dissipative stage are well approximated by linearized equations, 
\begin{equation}
\label{eq34}
\begin{cases}
\Delta u=p_2(\gamma-u_c-\alpha v_c) \\
\Delta v=p_2(1-v_c-\alpha w_c) \\
\Delta w=p_2(1-w_c-\alpha u_c)
\end{cases}
\end{equation}
These increments affect the quantity $A= {uvw}/{(u+v+w)^3}$. Assuming $u_c+v_c+w_c\approx1$, we derive the evolution of $\ln A$ over the dissipative stage with size $p_2$:

\begin{equation}
\label{eq35}
\Delta (\ln A) = p_2 [ (\gamma - 1) - (\alpha - 2) - 3(\gamma - 1) u_c + 3(\alpha - 2)(u_c v_c + v_c w_c + w_c u_c)].
\end{equation}
Expression (\ref{eq35}) is not sign conserving, but because the period $T_{\rm cons}$ of one conservative oscillation with fixed coexistence index $A$ is generically not commensurate with the duration $p_1$, one can expect that the points $(u_c,v_c,w_c)$ are uniformly distributed along the trajectory  with the fixed value of  $A$ (the change of $A$ gives the contribution of the higher order in $p_2$). Therefore, the total effect of dissipative stages can be estimated using the {\em average value} of $\Delta(\ln A)$ along the trajectory with the fixed $A$, i.e., as
\begin{equation}
\langle\Delta(\ln A)\rangle = p_2[(\gamma-1)-(\alpha-2)-3(\gamma-1)\langle u\rangle+3(\alpha-2)(\langle uv\rangle+\langle vw\rangle+\langle wu\rangle)].
\label{eq36}
\end{equation}
Here we use the notation
$$\langle f\rangle=\frac{1}{T_{cons}} \int_{0}^{T_{cons}}f(\tau)d\tau,$$
where $T_{cons}=T_{cons}(A)$ is the period of the oscillatory solution with a coexistence index $A$ at the conservative stage.

Now, since Eq.(\ref{eq3}) implies
\begin{eqnarray*}
d(\ln u)/{dt}&=&1-u-2v,\\
d(\ln v)/dt&=&1-v-2w,\\
d(\ln w)/{dt}&=&1-w-2u,
\end{eqnarray*}
after taking average over the period $T_{\rm cons}$ from both sides we obtain
\begin{equation}
\label{eq37}
\langle u \rangle=\langle v \rangle=\langle w \rangle=\frac{1}{3}.
\end{equation}
Since $v(t)=u(t-T/3),w(t)=u(t+T/3)$, it follows that 
$\langle u v \rangle = \langle v w \rangle = \langle w u \rangle$
so together with (\ref{eq37}):
\begin{equation}
\label{eq38}
\langle\Delta \ln A\rangle=cp_2 \  , \ c=(\alpha-2)(-1+9\langle u v \rangle)
\end{equation}
Generally, $c$ depends on $A$, and it can be seen that the sign of $c$ determines the stability of the heteroclinic cycle. Now, note that the action 
$$\langle f,g\rangle=\frac{1}{T_{\rm cons}}\int_0^{T_{\rm cons}}f(\tau)g(\tau)d\tau$$ is an inner product.
From Cauchy-Schwarz inequality, $|\langle uv\rangle| <\sqrt{\langle u^2\rangle}\sqrt{\langle v^2\rangle}$ ($u,v$ are linearly independent), and that $v(t)=u(t-T/3)$ one can show that:
\begin{equation}
\label{eq39}
|\langle u v\rangle|<|\langle u^2\rangle|.
\end{equation}

\begin{figure}[htbp]
\centering
\includegraphics[width=0.8\linewidth]{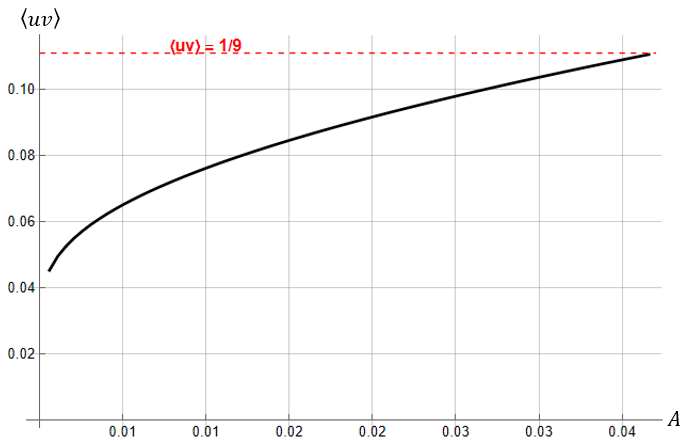}
\caption{\label{uv}\textbf{Mean product $\langle uv\rangle$ as a function of the coexistence index} $A$.
Trend of $\langle uv\rangle$ versus $A$ consistent with the sign of the averaged update for $\ln A$, supporting attraction toward the heteroclinic cycle in the appropriate regime.}
\end{figure}

From (\ref{eq11a}), $\dot{u}=u-u^2-2uv$, so after taking average on both sides, from (\ref{eq37}), we found that $\frac{1}{3}=\langle u^2 \rangle+2\langle uv \rangle$. Using (\ref{eq39}) we find that $\langle uv \rangle <1/9$ (see also Fig.\ref{uv} for numerical qualification). Therefore
\begin{equation}
\label{eq40}
 c<0
\end{equation}
and:
\begin{equation}
\label{eq41}
\frac{\langle A_{n+1}\rangle}{\langle A_n \rangle} \approx 1+cp_2,
\end{equation}
since $p_2\ll1$.
This implies that $A_n\rightarrow0$ as $n \to \infty$ which means that every initial condition that is not the coexistence state, is under sufficiently small $p_2$ eventually attracted by the heteroclinic cycle.

\begin{figure}[htbp]
\centering
\includegraphics[width=0.9\linewidth]{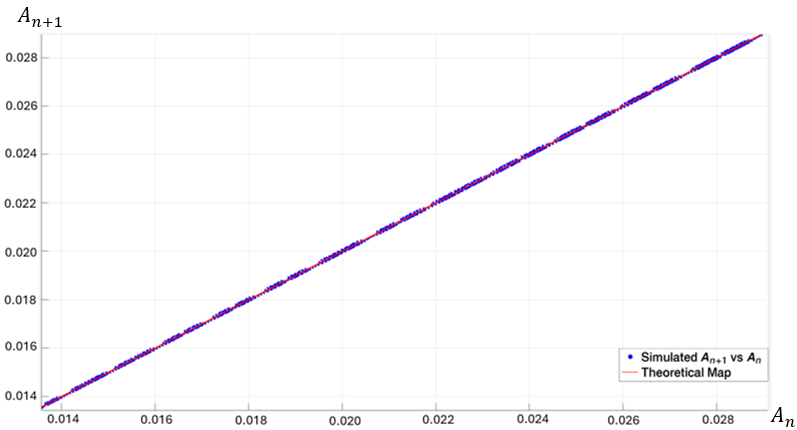}
\caption{\label{figure10}\textbf{Return plot of the discrete coexistence index} $A_{n+1}$ vs.\ $A_n$: simulation versus theory.}
Blue dots (simulation) and a red curve (averaged theory) are compared.
\end{figure}

\begin{figure}[htbp]
\centering
\includegraphics[width=0.8\linewidth]{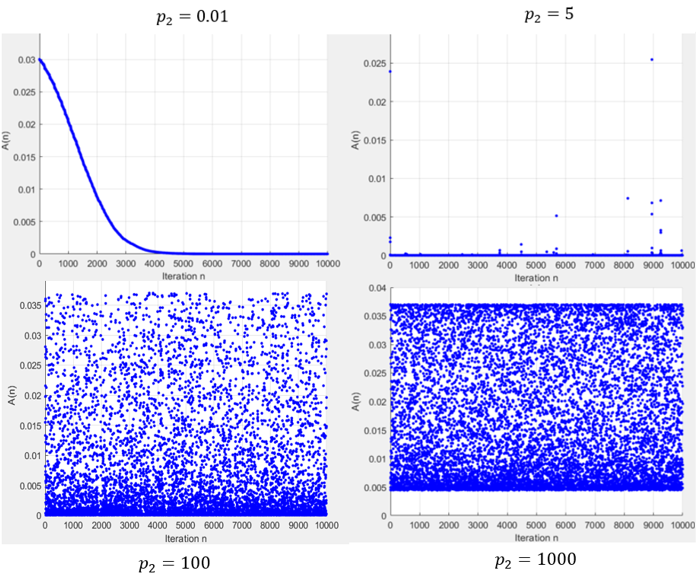}
\caption{\label{figure11}
\textbf{“Popcorn” dynamics of the  coexistence index} $A_n$ for varying dissipative season $p_2$ at fixed conservative season $p_1=1000$.}
As $p_2$ increases, rare upward bursts in $A_n$ become frequent and then irregular, signaling a transition from near-zero values to chaotic behavior.
\end{figure}

Note that the analysis above has been done for general initial conditions, that are not exactly on the heteroclinic cycle. The only exception is the initial condition $u=v=w=1/3$. This point is the coexistence state of the conservative stage, and not a periodic solution with nonzero $T_{\rm cons}$. A coexistence point of the conservative stage is not an equilibrium state during the dissipative stage. Therefore, even if the system starts at the coexistence point, in the dissipative season it is forced to leave this state, and the same mechanism applies in subsequent cycles. Thus, in the regime $p_2 \ll 1$, $p_1 \gg 1$, the heteroclinic cycle of the conservative dynamics becomes globally attracting. \\
Furthermore, as shown in \cite{Sorin2025}, for states near $A = 1/27$, we can expand $\langle u v \rangle \approx 1/9 - r$, with $r=\frac{1}{27}-A \ll 1$, leading to the approximation: $\langle u v\rangle=\frac{2}{27}+A$. Substituting it in (\ref{eq38}) gives that $c=(\alpha-2)(9\langle A_n\rangle-1/3)$, so that
\begin{equation}
\label{eq42}
\langle A_{n+1}\rangle =\langle A_n\rangle(1+p_2\frac{\alpha-2}{3}(27\langle A_n\rangle-1)).
\end{equation}
 This matches numerical simulations even beyond the vicinity of the coexistence state (see Fig.~\ref{figure10}).
 \\
 \\
As $p_2$ increases, the system changes its behavior. While for most $n$ $A_n$ is close to zero, occasional jumps in $A_n$ become more frequent, a phenomenon we call a "popcorn" (see Fig.~\ref{figure11}).

For sufficiently large $p_2$, the sequence ${A_n}$ appears chaotic, demonstrating a type of behavior not found in classical autonomous Lotka–Volterra or May–Leonard systems.

\section{Discussion}

\textbf{(What we did)}
\newline
We have explored a class of seasonally forced competitive Lotka–Volterra models, where the coefficients at the interaction terms are modeled as piecewise constant functions of time. 

\textbf{(What we achieved)}\newline
The presented approach explains how seasonal switching generates regimes that the autonomous models cannot: phase locking (synchronization windows), period-doubling cascades, irregular “popcorn” excursions of $A_n$, and robust chaos for sufficiently long dissipative stages. In the opposite limit (long conservative stage, short dissipative stage) drives $A_n$ to zero and makes the heteroclinic cycle globally attracting.  

\textbf{(Why this is relevant)} \newline
These findings indicate that seasonality, even when introduced as a relatively simple modulation, can fundamentally reshape ecological interactions and long-term evolutionary dynamics. 

\textbf{(Limitations and Extensions)}\newline
While our analysis focused on specific model families, the methods and insights developed here are applicable to a wider class of systems, including those with spatial structure or stochastic effects.
Nevertheless, we did not attempt data-based parameter identification. These are natural next steps rather than conceptual obstacles.

\textbf{(Compare with other works for benchmarking)}
\newline
The presented approach differs from classical “seasonal succession” studies in LV/May--Leonard models, which typically alternate growth and decay and focus on equilibria or heteroclinic cycles under decay phases. Here, both seasons are fully dynamical and the alternation is conservative$\leftrightarrow$dissipative, producing new regimes (locking, cascades, chaos) and a useful stroboscopic reduction. 

\textbf{(Perspectives)}\newline
Our future research will aim to characterize the structure and stability of the nontrivial solutions identified in this study, and to generalize the framework to more complex dynamical settings, such as systems with additional species, higher-dimensional phase spaces, or alternative forms of temporal forcing.

\FloatBarrier
\section{Conclusion}
\textbf{here please add shortest BOTTOMLINE TAKE-HOME MESSAGE}\newline
 Seasonal forcing of cyclic competition produces regimes that autonomous models cannot: synchronization windows, cascades to chaos, and when seasons are strongly imbalanced, global attraction to the heteroclinic cycle.
 
\textbf{(what was investigated?)}\newline
We studied a three-species system with periodic switching between a conservative stage and a dissipative stage, and analyzed the dynamics via a stroboscopic map and a coexistence index invariant on the May--Leonard manifold.

\textbf{(what was found)}\newline
 We found that season lengths select the long-term behavior: short dissipative seasons contract the invariant on average, driving the system toward the heteroclinic cycle; longer dissipative seasons introduce sensitive dependence, period-doubling, and robust chaos.
 
\textbf{(why relevant / important ; who should be interested})\newline
The presented approach shows how realistic season lengths can stabilize coexistence, create multiyear rhythms, or trigger chaos in rock--paper--scissors  population dynamics\index{Rock-Paper-Scissors(RPS)}. It offers testable predictions for phase locking and irregular fluctuations, and a compact toolset for nonautonomous ecological dynamics.

\end{document}